\documentclass[12pt,leqno]{article}
\title{Singular  integers and $p$-class group of cyclotomic fields    }
\author{Roland Qu\^eme}
\usepackage{amsmath,amsthm,amstext,amsbsy}
\newtheorem{thm}{Theorem}[section]

\newtheorem{lem}[thm]{Lemma}
\font\mathbb=msbm10
\newcommand{\N}{\mbox{\mathbb N}}

\newcommand{\Q}{\mbox{\mathbb Q}}
\newcommand{\Z}{\mbox{\mathbb Z}}

\newcommand{\modu}{\ \mbox{mod}\ }
\newcommand{\be}{\begin{equation}}
\newcommand{\ee}{\end{equation}}
\newcommand{\bd}{\begin{displaymath}}
\newcommand{\ed}{\end{displaymath}}
\newcommand{\bn}{\begin{enumerate}}
\newcommand{\en}{\end{enumerate}}

\usepackage{eucal}
\usepackage{eufrak}
\newcommand{\mk}{\mathfrak}

\date{2009 Oct 19}
\begin{document}
\maketitle
\abstract
Let $p$ be an irregular prime. Let $K=\Q(\zeta)$ be the $p$-cyclotomic field.
From Kummer and class field theory, there exist Galois extensions $S/\Q$ of degree $p(p-1)$ such that $S/K$ is a cyclic unramified extension of degree $[S:K]=p$.
We give an algebraic construction of the subfields $M$ of $S$ with degree $[M:\Q]=p$
and an explicit formula for the prime decomposition and ramification of the prime number $p$ in the extensions  $S/K$,  $M/\Q$ and $S/M$.
In the last section, we examine the consequences of these results for  the  Vandiver's conjecture.
This article is at elementary level on Classical Algebraic Number Theory.
%
%
\tableofcontents
%
\section{Some definitions}\label{s601191}
In this section   we give some  definitions and  notations on cyclotomic
fields,  $p$-class group, singular  numbers,  primary and non-primary, used in this paper.
\begin{enumerate}
\item
Let $p$ be an odd prime. Let $\zeta$ be a root of the polynomial equation $X^{p-1}+X^{p-2}+\dots+X+1=0$.
Let $K$ be the $p$-cyclotomic field $K=\Q(\zeta)$ and  $O_K$ its ring of integers.
Let $K^+$ be the maximal totally real subfield of $K$,
$O_{K^+}$  its ring of integers and  $O_{K^+}^*$ the  group of unit of $O_{K^+}$.
Let $v$ be a primitive root $\modu p$ and $\sigma: \zeta\rightarrow \zeta^v$ be a $\Q$-automorphism of $K$.
Let $G$ be the Galois group of the extension $K/\Q$.
Let ${\bf F}_p$ be the finite field of cardinal $p$ and   ${\bf F}_p^*$ its multiplicative group.
Let $\lambda=\zeta-1$. The prime ideal of $K$ lying over $p$ is $\pi=\lambda O_K$.
\item
Let $C_p$ be the $p$-class group of $K_p$ (the set of classes whose order is $1$ or $p$).
Let $r$ be the rank of $C_p$ seen as a ${\bf F}_p[G]$-module. If $r>0$ then $p$ is irregular.
Let $C_p^+$ be the $p$-class group of $K_p^+$. Let $r^+$ be the rank of $K_p^+$.  Then $C_p=C_p^+\oplus C_p^-$ where $C_p^-$
is the relative $p$-class group.
\item
$C_p$ is the direct sum of $r$ subgroups  $\Gamma_i$  of order $p$, each $\Gamma_i$
annihilated by a polynomial $\sigma-\mu_i\in {\bf F}_p[G]$ with
$\mu_i\in{\bf F}_p^*$,
\be\label{e611011}
C_p=\oplus_{i=1}^r \Gamma_i.
\ee
Then $\mu\equiv v^{n}\modu p$ with a natural integer
$n,\quad 1\leq n\leq p-2$.
\item
An  integer  $A\in O_K$ is said singular if $A^{1/p}\not\in K$
and if  there exists an   ideal $\mk a$  of $O_K$ such that
$A O_K =\mk  a^p$. Observe that, with this definition, a unit $\eta\in O_{K^+}^*$
with $\eta^{1/p}\not \in O_{K^+}^*$ is singular.
\item
A number $A\in K$ is said {\it semi-primary} if $v_\pi(A)=0$ and if there exists a natural integer $a$ such that
$A\equiv a\modu\pi^2$.
A number $A\in K$ is said {\it primary} if $v_\pi(A)=0$ and if there exists a natural integer $a$ such that
$A\equiv a^p\modu\pi^p$. Clearly a primary number is semi-primary.
A number $A\in K$ is said {\it hyper-primary} if $v_\pi(A)=0$ and if there exists a natural integer $a$ such that
$A\equiv a^p\modu\pi^{p+1}$. \en
%
\section{Some preliminary results}\label{s909121}
In this section we recall some properties of singular numbers given in Qu\^eme \cite{que}
in theorems  2.4 p. 4, 2.7 p. 7 and 3.1 p. 9. Let $\Gamma$ be one of the $r$ subgroups $\Gamma_i$ defined in relation
(\ref{e611011}).
\begin{enumerate}
\item\label{i609092}
\underline {If  $r^->0$ and  $\Gamma\subset  C_p^-$: }
then there exist  singular semi-primary integers $A$  with $A  O_K =\mk a^p$
where $\mk a$ is a {\bf non}-principal  ideal of $O_K$ and verifying   simultaneously
\begin{equation}\label{e512101}
\begin{split}
& Cl(\mk a)\in \Gamma,\ Cl(\mk a^{\sigma-\mu})=1,\\
& \sigma(A)=A^\mu\times\alpha^p,\quad \mu\in {\bf F}_p^*,\quad \alpha\in K,\\
&\mu\equiv v^{2m+1}\modu p, \quad m\in\N, \quad 1\leq m\leq \frac{p-3}{2},\\
&\pi^{2m+1} \ |\ A-a^p,\quad a\in\N,\quad 1\leq a\leq p-1.\\
\end{split}
\end{equation}
In that case we say that $A$ is a {\it negative} singular integer to point out that $Cl(\mk  a)\in C_p^-$.
Moreover, this number $A$ verifies
\begin{equation}\label{e512103}
A\times\overline{A}=D^p,
\end{equation}
for some integer $D\in O_{K^+}^*$.
\begin{enumerate}
\item
Either  $A$ is singular non-primary  with  $\pi^{2m+1} \ \|\ A-a^p$.
\item
Or  $A$ is singular  primary  with   $\pi^{p}\   |\ A-a^p$. In that case we know from class field theory  that $r^+>0$.
\end{enumerate}
(see Qu\^eme \cite{que} theorem 2.4 p. 4).
The singular primary negative numbers are interesting because
they  exist if and only if $h^+\equiv 0\modu p$ (the Vandiver conjecture is false).
\item\label{i609091}
\underline {If $r^+>0$ and    $\Gamma\subset C_p^+$: }
then there exist  singular semi-primary  integers $A$  with $A O_K =\mk a^p$
where $\mk a$ is a {\bf non}-principal  ideal of $O_K$ and verifying   simultaneously
\begin{equation}\label{e6012210}
\begin{split}
& Cl(\mk a)\in \Gamma,\ Cl(\mk a^{\sigma-\mu})=1,\\
& \sigma(A)=A^\mu\times\alpha^p,\quad \mu\in {\bf F}_p^*,\quad \alpha\in K,\\
&\mu\equiv v^{2m}\modu p, \quad m\in\Z, \quad 1\leq m\leq \frac{p-3}{2},\\
&\pi^{2m} \ |\ A-a^p,\quad a\in\Z,\quad 1\leq a\leq p-1,\\
\end{split}
\end{equation}
In that case we say that $A$ is a {\it positive } singular integer to point out that $Cl(\mk a)\in C_p^+$.
Moreover, this integer $A$ verifies
\begin{equation}\label{e512103}
\frac{A}{\overline{A}}=D^p,
\end{equation}
for some number  $D\in K_p^+$. If $h^+\equiv 0\modu p$ then $D\not=1$ is possible, for instance with
$\mk a=\mk q$ where $\mk q$ is a prime ideal of $O_K$,  $Cl(\mk q)\in C_p^+$ and $q\equiv 1\modu p$.
\begin{enumerate}
\item
Either  $A$ is singular non-primary with  $\pi^{2m} \ \|\ A-a^p$ .
\item
Or  $A$ is singular  primary with  $\pi^{p}\  |\ A-a^p$.
\end{enumerate}
(see Qu\^eme, \cite{que} theorem 2.7 p. 7).
\item\label{i609111}
\underline{If $\mu\equiv v^{2m}\modu p$ with $1\leq m\leq\frac{p-3}{2}$}:
then there exist singular units $A\in O_{K^+}^*$ with
\begin{equation}\label{e609091}
\begin{split}
& \sigma(A)=A^\mu\times \alpha^p,\quad \mu\in {\bf F}_p^*,\quad \alpha\in O_{K^+}^*,\\
&\mu\equiv v^{2m}\modu p, \quad m\in\Z, \quad 1\leq m\leq \frac{p-3}{2},\\
&\pi^{2m} \ |\ A-a^p,\quad a\in\Z,\quad 1\leq a\leq p-1,\\
\end{split}
\end{equation}
\begin{enumerate}
\item
Either $A$ is  non-primary with  $\pi^{2m} \ \|\ A-a^p$ .
\item
Or    $A$ is  primary with  $\pi^{p}\  |\ A-a^p$.
\en
(see Qu\^eme, \cite{que} theorem 3.1 p. 9).
\en
%
The sections \ref{s909271}, \ref{s909273}  and \ref{s909272} are,  for a large part,   a reformulation of Hilbert theory of {\it  Kummer Fields}, see \cite{hil} paragraph 125 p. 225.
\section{Singular $K$-extensions}\label{s909271}
\paragraph{Some Definitions}:
\bn
\item
In this section, let us denote  $\Gamma$  one of the $r$ subgroups of order $p$ of $C_p$ defined by relation (\ref{e611011}).
Let $A$ be a singular semi-primary  integer,    negative or  positive,  verifying   respectively  the relations
(\ref{e512101}) or (\ref{e6012210}) .
We call $S =K(A^{1/p})/K$ a  singular negative, respectively  positive  $K$-extension if $\Gamma\in C_p^-$, respectively $\Gamma\in C_p^+$.
\item
Let $A$ be a singular unit verifying the relation  (\ref{e609091}).
We call $S =K(A^{1/p})/K$ a  singular unit  $K$-extension.
\item
A singular $K$-extension $S=K(A^{1/p})$ is said primary or  non-primary
if the singular number $A$  is primary or non-primary.
\item
If $S$ is  primary then the extension  $S/K$ is, from Hilbert class field theory,
the cyclic unramified extension of degree $p$ corresponding to $\Gamma$.
\item Observe that the extensions $S/\Q$  are Galois extensions of degree $p(p-1)$.
\en
%
\begin{lem}\label{l609111}
There is one and only one singular negative $K$-extension  corresponding to a group $\Gamma\subset C_p^-$.
\begin{proof}
For $\Gamma$ given let us consider two singular negative $K$-extensions $S/K$ and $S^\prime/K$.
$A O_K=\mk a^p$ and $A^\prime O_K= \mk a^{\prime p}$. The polynomial  $\sigma-\mu$ annihilates
$<Cl(\mk a>$ and $<Cl(\mk a^\prime>$.
Then $<Cl(\mk a)>\ =\ < Cl(\mk a^\prime)>\ =\Gamma$,
thus there exists $n,\ 1\leq n\leq p-1$ such that $Cl(\mk a^n)=Cl(\mk a^\prime)$.
Therefore $A^n= A^\prime\times\gamma^p\times \varepsilon,\ \varepsilon\in O_{K}^*,\ \gamma\in K$.
It follows, from $A\overline{A} =D^p$ and $A^\prime \overline{A}^\prime=D^{\prime p}$ with $D, D^\prime\in O_{K^+}$,
that $\varepsilon\overline{\varepsilon} \in O_{K^+}^{* p}$.
Therefore $\varepsilon= \zeta^w\varepsilon_1^p,\ \varepsilon_1\in O_{K^+}^*$.
Then $A^n=A^\prime\gamma^p\zeta^w\varepsilon_1^p\zeta^w$.
$A$ and $A^\prime$ are semi-primary, thus  it follows that $w=0$. Therefore $K(A^{1/p})= K(A^{\prime 1/p})$.
\end{proof}
\end{lem}
\paragraph{Remark:} Observe that we consider in this article only  singular semi-primary numbers. Let $A$ be a singular semi-primary number. Then
$A^\prime= A\zeta$ is not semi-primary and $K(A^{1/p})\not=K(A^{1/p})$.
%
%
\begin{lem}\label{l609113}
If $\mu\not=1$ and $\mu^{(p-1)/2}\equiv 1\modu p$   there is one and only one  singular  unit $K$-extension  $S/K$ depending only on $\mu$.
\begin{proof}
The subgroup of $O_{K^+}^*/O_{K^+}^{* p}$   annihilated by $\sigma-\mu$ is of order $p$ and the rank of $O_{K^+}^*/O_{K^+}^{* p}$
is  $\frac{p-3}{2}$.
\end{proof}
\end{lem}
%
%
\begin{lem}\label{l609091}
$\pi$ is the only prime which can ramify in the singular $K$-extension $S/K$ and the relative discriminant of $S/K$ is a power of $\pi$.
\begin{proof}
$S/K$ is unramified except possibly at $\pi$, (see for instance  Washington \cite{was} exercise 9.1 (b) p. 182).
The result for relative discriminant follows.
\end{proof}
\end{lem}
%
\begin{lem}\label{l609092}$ $
\bn
\item
There are $r^+$ singular primary negative   $K$-extensions $S/K$.
\item
There are $r^--r^+$ singular non-primary  negative  $K$-extensions $S/K$.
\en
\begin{proof}
The first part results of classical theory of $p$-Hilbert class field applied to the field $K$ and of previous definition of singular $K$-extensions $S_\mu$
(see for instance the result of Furtwangler in Ribenboim \cite{rib} (6C) p. 182) and the second part is an immediate consequence of the first part.
\end{proof}
\end{lem}
%
\section{Singular $\Q$-fields }\label{s909273}
Let $A$ be a semi-primary  integer, negative (see definition (\ref{e512101})), positive (see definition (\ref{e6012210}))
or unit (see definition (\ref{e609091})).
Let $\omega$ be an algebraic number defined by
\be\label{e610096}
\omega= A^{(p-1)/p}.
\ee
We had chosen this definition instead of $\omega=A^{1/p}$ because $A^{p-1}\equiv 1\modu \pi$ simplifies computations.
Then $S=K(\omega)$ is  the corresponding singular $K$-extension.
Observe that this definition implies that
$\omega\in O_{S}$ ring of integers of $S$.
%
%
\begin{lem}\label{l609101}
Suppose that $S/K$ is a singular primary $K$-extension.
Let $\theta :\omega\rightarrow\omega\zeta$ be a $K$-isomorphism of the field $S$.
Then, $A$ is hyperprimary and there are $p$ prime ideals  of $O_{S}$ lying over $\pi$.
There  exists a prime ideal $\pi_0$ of $O_{S}$ lying over $\pi$ such that
the $p$ prime ideals $\pi_n=\theta^n(\pi_0),\ n=0,\dots,p-1$ of $O_{S}$ lying over $\pi$ verify the congruences
\be\label{e609043}
 \begin{split}
 &\pi_0^ 2\ |\ \omega-1,\\
 &\pi_n\ \|\ \omega-1, \ \dots, n=1,\dots,p-1.\\
\end{split}
 \ee
  \begin{proof}$ $
\bn
\item
From Hilbert class field theory and Principal Ideal Theorem the prime principal ideal $\pi$ of $K$ splits totally in the extension $S/K$.
The ideal $\pi$ does not correspond to the case III.c in Ribenboim \cite{rib} p. 168 because $\pi$ is not ramified in $S/K$.
The ideal $\pi$ does not correspond to the case III.b in Ribenboim \cite{rib} p. 168 because $\pi$ is not inert  in $S/K$.
Therefore $\pi$ corresponds to the case III.a and it follows that there exists $a_1\in O_K$ such that $A\equiv a_1^p \modu \pi^{p+1}$.
Therefore there exists $a\in\Z$  such that $a\equiv a_1\modu\pi$ and  $A\equiv a^p \modu \pi^{p+1}$, thus $A$ is a singular hyper-primary number
and $A^{p-1}\equiv 1\modu\pi^{p+1}$.
\item
Then $\omega^p-1\equiv 0\modu\pi^{p+1}$.
Let $\theta :\omega\rightarrow \omega\zeta$ be a $K$-automorphism of the field  $S$.
Let $\Pi^\prime$ is any  of the $p$ prime ideals of
$O_{S}$ lying over $\pi$.
Then  $\pi O_{S}=\prod_{n^\prime=0}^{p-1}\pi_{n^\prime}$ where
$\pi_n^\prime=\theta^n(\Pi^\prime),\ n=0,\dots,p-1$ are the $p$ prime ideals of $O_{S}$ lying over $\pi$.
\item
From $A^{p-1}\equiv 1\modu \pi^{p+1}$ we see that
\bd
\omega^p-1=\prod_{n=0}^{p-1}(\omega\zeta^{-n}-1)\equiv 0\modu\pi_0^{\prime p+1}\pi_1^{\prime p+1}\dots\pi_{p-1}^{\prime p+1}.
\ed
It follows   that there exists a prime ideal $\Pi$ of $O_{S}$ lying over $\pi$
such that $\omega-1\equiv 0\modu \Pi^{ 2}$ because there exists $l$ such that $\omega\zeta^l-1\equiv 0\modu \pi_0^{\prime 2})$
so $\Pi=\theta^{-l}(\pi_0^\prime$,
and that  $\Pi\ \|\ \omega\zeta^{-n}-1$ for $n=1,\dots,p-1$ because
$\Pi\ \|\ \zeta-1$.
Let us note $\pi_n=\theta^n(\Pi)$ for $n=0,\dots,p-1$.
It follows that, for $n=1,\dots, p-1$,
 $\pi_n\ \|\ \omega\zeta^{n}\zeta^{-n}-1 $ and so
 \be\label{e609043}
 \begin{split}
 & \pi_0^2\ |\ \omega-1,\\
 &\pi_n\ \|\ \omega-1, \ \dots, n=1,\dots,p-1.\\
\end{split}
 \ee
 \en
 \end{proof}
 \end{lem}
%
%
\begin{lem}\label{l609121}
Suppose that $S/K$ is a singular non primary  $K$-extension.
Let $\Pi$ be the prime of $S$ lying over $\pi$. Then $\Pi\ |\ \omega-1$.
\begin{proof}$ $
The extension $S/K$ is ramified therefore $\pi O_{S}=\Pi^p$.
$A^{p-1}\equiv 1\modu\pi^n$ for some $n>1$ and so $\omega^p-1\equiv 1\modu\Pi^{np}$ because $\pi O_{S}=\Pi^p$.
Therefore $\omega\equiv 1\modu\Pi$.
\end{proof}
\end{lem}
%
We know that there are $p$ different automorphisms   of the field $S$ extending the $\Q$-automorphism $\sigma$ of the field $K$.
\begin{lem}\label{l609095}
There exists  an automorphism  $\sigma_\mu$ of  $S/\Q$ extending $\sigma$  such that
\be\label{e609095}
\omega^{\sigma_\mu-\mu}\equiv 1\modu \pi^2.
\ee
\begin{proof}
\item
From $\sigma(A)=A^\mu\alpha^p$ there exist $p$ different automorphisms $\sigma_{(w)},\ w=0,\dots,p-1,$ of the field  $S$
extending the $\Q$-isomorphism
$\sigma$ of the field $K$,  defined by
\be\label{e610137}
\sigma_{ (w)}(\omega)=\omega^\mu\alpha^{p-1}\zeta^w,
\ee
for  natural numbers $w=0,1,\dots,p-1$.
There exists one and only one  $w$ such that $\alpha^{p-1}\times \zeta^w$ is a semi-primary number (or
$\alpha^{p-1}\times\zeta^w\equiv 1\modu\pi^2$). Let us set $\sigma_\mu=\sigma_{(w)}$ to emphasize the role of $\mu$
Therefore we get
\be\label{e609041}
\sigma_\mu(\omega)\equiv \omega^\mu\modu\pi^2,
\ee
because $\omega,\sigma_\mu(\omega)\in O_{S}$.
\end{proof}
\end{lem}
%
\begin{lem}\label{l609096}
$\sigma_\mu^{p-1}(\omega)=\omega$.
\begin{proof}
We have $\sigma_\mu^{p-1}(A)= \sigma^{p-1}(A)= A$ therefore there exists a natural integer $w_1$ such that
$\sigma_\mu^{p-1}(\omega)=\omega\times \zeta^{w_1}$.
We have proved in relation (\ref{e609041}) that
\be\label{e609103}
\sigma_\mu(\omega)\equiv \omega^\mu \modu \pi^2,
\ee
thus $\sigma_\mu^{p-1}(\omega)\equiv\omega^{\mu^{p-1}}\equiv\omega\times  A^{(p-1)(\mu^{p-1}-1)/p}\equiv \omega\modu \pi^2$
which implies that $w_1=0$ and that $\sigma_\mu^{p-1}(\omega)=\omega$.
\end{proof}
\end{lem}
%
Let us define $\Omega\in O_{S}$ ring of integers of $S$ by the relation
\be\label{e609061}
\Omega=\sum_{i=0}^{p-2}\sigma_\mu^i(\omega).
\ee
%
\begin{thm}\label{l609097}
$M=\Q(\Omega)$ is a field with $[M:\Q]= p$, $[S: M]=p-1$  and $\sigma_\mu(\Omega)=\Omega$.
\begin{proof}$ $
\bn
\item
Show that $\Omega\not=0$: If $S/K$ is unramified, then $\omega\equiv 1\modu\pi$ implies with definition of $\Omega$ that
$\Omega\equiv p-1\modu \pi$ and so $\Omega\not=0$. If $S/K$ is ramified, then $\omega\equiv 1\modu\Pi$
implies with definition of $\Omega$ that
$\Omega\equiv p-1\modu \Pi$ because $\sigma_\mu(\Pi)=\Pi$ and so $\Omega\not=0$.
\item
Show that $\Omega\not\in K$: from $\sigma_\mu(\omega)=\omega^\mu\alpha^{p-1}\zeta^w$ we get
\bd
\Omega=\sum_{i=0}^{p-2}\omega^{\mu^i\mod p}\times \beta_i,
\ed
with $\beta_i\in K$.
Putting together terms of same degree   we get  $\Omega=\sum_{j=1}^{p-1}\gamma_j\omega^j$ where $\gamma_j\in K$ are not all null
because $\Omega\not=0$.
$\Omega\in K$
should imply the polynomial equation $\sum_{j=1}^{p-1}\omega^{j}\times\gamma_j-\gamma=0$ with $\gamma\in K$,  not possible
because the minimal polynomial equation of $\omega$ with coefficients in $K$ is $\omega^p-A^{p-1}=0$.
\item
Show that $M=\Q(\Omega)$ verifies $M\subset S$ with $[M:\Q]=p$ and  $[S:M]=p-1$:
$S/\Q$ is a Galois extension with $[S:\Q]=(p-1)p$.
Let $G_S$ be the Galois group of $S/\Q$. Let $<\sigma_\mu>$ be the subgroup of $G_S$
generated by the automorphism $\sigma_\mu\in G_S$. We have seen in lemma \ref{l609096} that $\sigma_\mu^{p-1}(\omega)=\omega$.
In the other hand $\sigma_\mu^{p-1}(\zeta)=\zeta$ and $\sigma_\mu^n(\zeta)\not=\zeta$
for $n<p-1$ and so $<\sigma_\mu>$ is of order $p-1$.
\item
From fundamental theorem of Galois theory,  there is a fixed field $M=S^{<\sigma_\mu>}$  with
$[M:\Q]= [G_S:<\sigma_\mu>]=p$.
From $\sigma_\mu(\Omega)=\Omega$ seen and from definition relation (\ref{e609061}) it follows that $\Omega\in M$ and from $\Omega \not\in K$
it follows that $M=\Q(\Omega)$.
Thus  $S=M(\zeta)$
and    $\omega\in S$ can be written
\be\label{e608093}
\omega=1+\sum_{i=0}^{p-2} \omega_i\lambda^i,\ \omega_i\in M.
\ee
with $\lambda=\zeta-1$ and with  $\sigma_\mu(\omega_i)=\omega_i$  because $\sigma_\mu(\Omega)=\Omega$.
\en
\end{proof}
\end{thm}
%
\paragraph{Some  definitions:}
The field  $M\subset S$ is called a singular $\Q$-field. In the sequel of this paper we are studying some algebraic properties and
ramification of singular $\Q$-fields $M$. A singular $\Q$-field $M$ is said primary (respectively  non-primary)
if $S$ is a singular primary (respectively  non-primary) $K$-extension.
%
\section{Algebraic properties  of singular $\Q$-fields}\label{s909272} 
\bn
\item
From Galois theory  there are $p$ subfields $M_i,\ i=0,\dots,p-1,$  of $S$  of degree $[M_i:\Q]=p$.
\item
The extension $S/\Q$ is Galois.
Let   $\theta :\omega\rightarrow\omega\zeta$ be a $K$-automorphism of $S$.
There are $p$ automorphisms $\sigma_i,\ i=0,\dots,p-1,$ of $S$
extending the $\Q$-automorphism $\sigma$ of $K$
verifying  $\sigma_i(\theta^i(\omega))=(\theta^i(\omega))^\mu\beta$
for the  semi-primary $\beta\in K$.
\item
We have defined in relation (\ref{e609061}) $\Omega=\sum_{k=0}^{p-2} \sigma_\mu^k(\omega)$.
For $i=1,\dots,p-1$ we can define similarly
$\Omega_i=\sum_{k=0}^{p-2} \sigma_\mu^k(\theta^i(\omega))$.
Then we  show in following result  that the  fields $M_i$ can be explicitly defined by $M_i=\Q(\Omega_i),\ i=0,\dots,p-1$.
\en
%
\begin{lem}\label{l6090981}
The singular $\Q$-fields $M_i=\Q(\Omega_i),\dots,  p-1$, are the $p$ subfields of degree $p$  of the singular $K$-extension $S/K$.
\begin{proof}$ $
\bn
\item
We set  here $\sigma_0 =\sigma_\mu$ and $M_0=M$.
Show that the fields $ M_0,M_1,\dots,M_{p-1}$ are pairwise different:
$\sigma_i(\theta^i(\omega))=(\theta^i(\omega))^\mu\beta$, hence
$\sigma_i(\omega\zeta^i)=(\omega\zeta^i)^\mu\beta$, hence
\begin{equation}\label{e611083}
\sigma_i(\omega)=\omega^\mu\beta\zeta^{i(\mu-v)}.
\end{equation}
Suppose that $M_i=M_{i^\prime}$: then the subgroups $<\sigma_i>$ and $<\sigma_{i^\prime}>$ of $Gal(S/\Q)$  corresponding to the fixed fields
$M_i$ and $M_{i^\prime}$ are equal. Therefore there exists a natural integer $l,\ 1\leq l\leq p-2$ coprime with $p-1$ such that
$\sigma_{i^\prime}=\sigma_i^l$.
\item
$\sigma_{i^\prime}(\zeta)=\sigma_i^l(\zeta)$, hence $\zeta^{v}= \zeta^{v^l}$, hence $v\equiv v^l\modu p$,
hence $v^{l-1}\equiv 1\modu p$, hence $l-1\equiv 0\modu p-1$
and therefore $l\equiv 1\modu p-1$. In the other hand $1\leq l\leq p-2$, thus $l=1$ and
$\sigma_i(\omega)=\sigma_{i^\prime}(\omega)$. From relation (\ref{e611083}) this implies that  $i=i^\prime$.
\en
\end{proof}
\end{lem}
%
In the following theorem, we give an explicit computation of $\Omega_i$ for $i=0,\dots,p-1$.
Let us denote  $\mu_k$ for $\mu^k\modu p$.
\begin{lem}\label{t610201}$ $
The subfields of degree $p$ of the singular $K$-extension $S$ are
the singular $\Q$-fields $M_i=\Q(\Omega_i),\ i=0,\dots,p-1,$
where
\be\label{e610201}
\begin{split}
&\Omega_i=\theta^i(\Omega)= \sum_{k=0}^{p-2}\omega^{\mu^k}\beta^{(\sigma^k-\mu^k)/(\sigma-\mu)}\zeta^{i\mu^k},\\
&\Omega_i=\theta^i(\Omega)= \sum_{k=0}^{p-2}\omega^{\mu_k}A^{(p-1)(\mu^k-\mu_k)/p}\beta^{(\sigma^k-\mu^k)/(\sigma-\mu)}\zeta^{i\mu^k}.\\
\end{split}
\ee
\begin{proof}
We start of
$\Omega_i=\sum_{k=0}^{p-2}\sigma_i^k(\theta^i(\omega))$
and we  compute $\sigma_i^k(\theta^i(\omega))$.
Let us note $\varpi_i=\theta^i(\omega)$.
$\sigma_i(\varpi_i)=\varpi_i^\mu\beta$, hence
$\sigma_i^2(\varpi_i)=\sigma(\varpi_i)^\mu\sigma(\beta)=(\varpi_i^\mu\beta)^\mu\sigma(\beta)=\varpi_i^{\mu^2}\beta^{\sigma+\mu}$.
Pursuing up to $k$,
we get $\sigma_i^k(\varpi_i)=\varpi_i^{\mu^k}\beta^{(\sigma^k-\mu^k)/(\sigma-\mu)}$.
But $\varpi_i^{\mu^k}=(\omega\zeta^i)^{\mu^k}=\omega^{\mu^k}\zeta^{i\mu^k}$.
We can also compute at first $\Omega=\sum_{k=0}^{p-2}\omega^{\mu^k}\beta^{(\sigma^k-\mu^k))/(\sigma-\mu)}$
and then verify directly that $\theta^i(\Omega)=\Omega_i$.
Then $\omega^{\mu^k}=\omega^{\mu_k} A^{(p-1)(\mu^k-\mu_k)/p}$.
\end{proof}
\end{lem}
%
%
\section{The ramification in the singular primary $\Q$-fields}
\bn
\item
Observe at first that the case of singular non-primary $\Q$-fields can easily be described. The extension $S/K$ is fully ramified at $\pi$,
so $p O_{S}=\pi_{S}^{p(p-1)}$. Therefore there is only one prime ideal $\mk p$ of $M$ ramified
with $p O_{M}=\mk p^p$.
\item
The end of this section deals with the ramification of singular primary  $\Q$-fields  $M$.
In that case  $S/K$ is a cyclic unramified extension and there are $p$ prime ideals in $S/K$ over $\pi$.
\en
%
\begin{lem}\label{l609102}
$\sigma_\mu(\pi_0)=\pi_0$
\begin{proof}
From relation (\ref{e609103}) $\sigma_\mu(\omega)\equiv\omega^\mu\modu\pi^2$. From lemma \ref{l609101} $\omega\equiv 1\modu\pi_0^2$
and so $\sigma_\mu(\omega)\equiv\omega^\mu\equiv 1\modu\pi_0^2$.  Then $\omega\equiv 1\modu \sigma_\mu^{-1}(\pi_0)^2$.
If $\sigma_\mu^{-1}(\pi_0)\not=\pi_0$ it follows that $\omega\equiv 1\modu\pi_0^2\times\sigma^{-1}(\pi_0^2)$, which contradicts lemma
\ref{l609101}.
\end{proof}
\end{lem}
%
\begin{lem}\label{l301281}
Let $\pi_k=\theta^k(\pi_0)$ for any $k\in\N,\quad 1\leq k\leq p-1$.
Then $\sigma_\mu(\pi_k)=\pi_{n_k}$ with $n_k\in\N, \quad n_k\equiv k\times v \mu^{-1}\modu p$.
\begin{proof}$ $
\bn
\item
From $\pi_0^{ 2}\  |\ (\omega-1)$, it follows that $\theta^k(\pi_0^{ 2})=\pi_k^{ 2}\ |\ (\omega\zeta^k-1)$.
Then
\bd
\sigma_\mu(\pi_k)^{ 2}\ |\ (\sigma_\mu(\omega)\times\zeta^{v k}-1).
\ed
\item
We have $\sigma_\mu(\pi_k )=\pi_{k+l_k}$ for some $l_k\in\N$ depending on $k$.
From relation (\ref{e609103}) we know that $\sigma_\mu(\omega)\equiv \omega^\mu \modu\pi^2$.
Therefore
\begin{displaymath}
\pi_{k+l_k}^{ 2}\ |\ (\omega^\mu\times\zeta_p^{v k}-1).
\end{displaymath}
\item
In an other part by the $K$-automorphism $\theta^{k+l_k}$ of $S$ we have
\bd
\pi_{k+l_k}^{ 2}\ |\ (\omega\times\zeta^{k+l_k}-1),
\ed
so
\begin{displaymath}
\pi_{k+l_k}^{ 2}\ |\ (\omega^\mu\times\zeta^{\mu(k+l_k)}-1).
\end{displaymath}
\item
Therefore
$\pi_{k+l_k}^{ 2}\ |\  \omega^\mu(\zeta^{v k}-\zeta^{\mu(k+l_k)})$,
and so
\begin{displaymath}
\pi_{k+l_k}^{ 2}\ |\ (\zeta^{v k}-\zeta^{\mu(k+l_k)}),
\end{displaymath}
\item
This implies that  $\mu(k+l_k)-v k\equiv 0\modu p$, so
$\mu l_k+k(\mu-v)\equiv 0\modu p$ and
finally that
\begin{displaymath}
l_k\equiv k\times\frac{v-\mu}{\mu},
\end{displaymath}
where we know that $v-\mu\not\equiv 0\modu p$ from Stickelberger relation. Then
$n_k\equiv k+k\times\frac{v-\mu}{\mu}=k\times\frac{v}{\mu}\modu p$, which achieves the proof.
\en
\end{proof}
\end{lem}
%
\begin{lem}\label{l609271}$ $
\bn
\item
If $S/K$ is a  singular primary negative extension then $\sigma_\mu^{(p-1)/2}(\pi_k)=\pi_k$.
\item
If $S/K$ is a singular primary positive or unit  extension then $\sigma_\mu^{(p-1)/2}(\pi_k)=\pi_{n-k}$.
\en
\begin{proof}
From lemma \ref{l301281} we have $\sigma_\mu^{(p-1)/2}(\pi_k)= \pi_{k^\prime}$ with
$k^\prime\equiv k v^{(p-1)/2}\mu^{-(p-1)/2}$.
 If $S/K$ is negative then  $v^{(p-1)/2}\mu^{-(p-1)/2}\equiv 1\modu p$ and if
$S/K$ is positive or unit then  $v^{(p-1)/2}\mu^{-(p-1)/2}\equiv -1\modu p$
and the result follows.
\end{proof}
\end{lem}
%
\begin{lem}\label{l609105}
The length of the orbit of the action of the group $<\sigma_\mu>$ on $\pi_0$ is $1$
and the length of the orbit of the action of the group $<\sigma_\mu>$ on
$\pi_i, \  i=1,\dots,p-1$ is $d$ where $d$ is the order of $v\mu^{-1}\modu p$.
\begin{proof}
For $\pi_0$ see lemma \ref{l609102}.
For $\pi_k$ see lemma \ref{l301281}: $\sigma_\mu(\pi_k)=\sigma(\pi_{n_k})$ with
$n_k\equiv v\mu^{-1}\modu p$, then $\sigma_\mu^2(\pi_k)=\sigma(\pi_{n_{k_2}})$ with
$n_{k_2}\equiv k v^2\mu^{-2}\modu p$ and finally $n_{k_d}\equiv k \modu p$.
\end{proof}
\end{lem}
%
The only prime ideals of $M/\Q$ ramified are lying over $p$.
The prime ideal of $K$ over p is $\pi$.
To avoid cumbersome notations,
the prime ideals of $S$ over $\pi$ are noted here
$\Pi$ or $\Pi_i=\theta^i(\Pi_0), \ i=1,\dots,p-1$,
and the prime ideals of $M$ over $p$ are noted $\mk p$ or  $\mk p_j,\ j=1\dots,\nu$ where   $\nu+1$ is the number of such ideals.
%
\begin{thm}\label{t609101}$ $
Let  $d$ be  the order of $v\mu^{-1}\modu p$.
There are $\frac{p-1}{d}+1$ prime ideals in the singular primary $\Q$-field $M$ lying over $p$.
Their  prime decomposition and ramification is:
\bn
\item
$ e(\mk p_0/p\Z) =1.$
\item
$e(\mk p_j/p\Z)=d$ for all $j=1,\dots,\frac{p-1}{d}$  with $d>1$.
\en
\begin{proof}$ $
\bn
\item
\underline{preparation of the proof}
\bn
\item
The inertial degrees verifies
$f(\pi/p\Z)=1$ and $f(\Pi/\pi)=1$ and so $f(\Pi/p\Z)=1$.
Therefore, from multiplicativity of degrees in extensions,
it follows that
$f(\mk p/p\Z)=f(\Pi/\mk p)=1$
where $\Pi$ is lying over $\mk p$.
\item
$e(\pi/p\Z)=p-1$ and $e(\Pi/\pi)=1$ and so $e(\Pi/p\Z)=p-1$.
\item
Classically, we get
\be\label{e609274}
\sum_{j=0}^\nu e(\mk p_j/p\Z)=p,
\ee
where $\nu+1$ is the number of prime ideals of $M$ lying over $p$ and
where $e(\mk p_j/p\Z)$ are ramification indices dividing $p-1$ because, from multiplicativity of degrees in extensions,
$e(\mk p_j/p\Z)\times e(\Pi/\mk p_j)=p-1$.
\en
\item
\underline{Proof}
\bn
\item
The extension $S/M$ is Galois of degree $p-1$, therefore
the number of prime ideals $\Pi$ lying over one $\mk p$ is $\frac{p-1}{e(\Pi/\mk p)}=e(\mk p/p\Z)$.
\item
Let $c(\Pi)$ be the orbit of $\Pi$ under the action of the group $<\sigma_\mu>$ of cardinal $p-1$ seen in lemma \ref{l609105}.
If $\Pi=\pi_0$ then the orbit $C_\Pi$ is of length $1$.
If $\Pi\not=\pi_0$ then the orbit $C_\Pi$ is of length $d$.
If $C_\Pi$ has one ideal lying over $\mk p$ then it has all its $d$ ideals lying over $\mk p$ because $\sigma_\mu(\mk p)=\mk p$.
This can be extended to   all $\Pi^\prime$ lying over $\mk p$ with $C_{\Pi^\prime}\not= C_\Pi$  and it follows that
when $\Pi\not=\pi_0$ then $d\ |\  e(\mk p/p\Z)$, number of ideals of $S$ lying over $\mk p$.
There is  one $\mk p$ with $e(\mk p/p\Z)=1$ because $C_{\pi_0}$ is the only orbit with one element.
\item
The extension $S/M$ is cyclic of degree $p-1$. There exists one field $N$ with
$M\subset N\subset S$ with degree  $[N:M]=\frac{p-1}{d}$.
If there were   at least two different prime  ideals $\mk p_{1}^\prime $ and $\mk p_{2}^\prime$ of $N$ lying over $\mk p$,
it should follow
that $\mk p_{2}^\prime=\sigma_\mu^{j d}(\mk p_{1}^\prime)$ for some $j,\ 1\leq j\leq d-1$
because the Galois group of $S/M$ is $<\sigma_\mu>$ and the Galois group of $N/M$ is $<\sigma_\mu^{d}>$.
But, if a prime ideal $\pi_k$ of $S$ lies over $\mk p_{1}^\prime$ then $\sigma_\mu^d(\pi_k)$ should lie over
$\mk p_{2}^\prime$.
From lemma \ref{l609105},  $ \sigma_\mu^{d}(\pi_k)=\pi_k$ should  imply that $\mk p_{2}^\prime=\mk p_{1}^\prime$, contradiction.
Therefore the only possibility is that $\mk p$ is fully ramified in $N/M$ and
thus  $\mk p O_{N}=\mk p^{\prime (p-1)/d}$.
Therefore $e(\mk p^\prime/\mk p)=\frac{p-1}{d}$ and so
$e(\mk p^\prime/p\Z)=e(\mk p/p\Z)\times \frac{p-1}{d}\ |\ p-1$ and thus $e(\mk p/p\Z)\ |\ d$.
From previous result it follows that $e(\mk p/p\Z)=d$. Then $d>1$ because $\mu-v\not\equiv 0\modu p$ from Stickelberger theorem.
There are $\frac{p-1}{d}+1$  prime ideals $\mk p_i$ because, from relation (\ref{e609274})
$p= 1+\sum_{i=1}^\nu e(\mk p_i/p\Z)=1+\nu \times d$.
\en
\en
\end{proof}
\end{thm}
%
\paragraph{Example:}
let us consider the case of prime numbers $p$ with  $\frac{p-1}{2}$  prime.
\bn
\item
\underline{Singular primary negative  $\Q$-fields}

Here $\mu^{(p-1)/2}\equiv -1\modu p$ and
$d\in\{2,\frac{p-1}{2},p-1\}$.
Straightforwardly  $d=2$ is not possible: $\mu^2\equiv v^2\modu p$, then $\mu+v\equiv 0\modu p$ because
$\mu\not\equiv v\modu p$, then $\mu^{(p-1)/2}+v^{(p-1)/2}\equiv 0\modu p$, contradiction
because $\mu^{(p-1)/2}=v^{(p-1)/2}=-1$. $d=p-1$ is not possible because  $\mu^{(p-1)/2}-v^{(p-1)/2}\equiv 0\modu p$.
Therefore $d= \frac{p-1}{2}$, so  the ramification of $p$ in the singular $\Q$-field $M$ is
$e(\mk p_0/p\Z)=1$ and $e(\mk p_1/p\Z)=e(\mk p_2/p\Z)=\frac{p-1}{2}$.
\item
\underline{Singular primary positive   $\Q$-extensions and primary unit $\Q$-fields}

Here $\mu^{(p-1)/2}\equiv 1\modu p$ and
$d\in\{2,\frac{p-1}{2},p-1\}$.
$d=2$ is not possible : $\mu^2-v^2\equiv 0\modu p$ then $\mu+v\equiv 0\modu p$ so $\mu\equiv v^{(p+1)/2}\modu p$,
so $B_{p-(p+1)/2}=B_{(p+1)/2}\equiv 0\modu p$
where $B_{(p+1)/2}$ is a Bernoulli Number, contradiction because $B_{(p+1)/2}\not\equiv 0\modu p$.
$d=\frac{p-1}{2}$ is not possible because $\mu^{(p-1)/2}\equiv 1\modu p$ and $v^{(p-1)/2}\equiv -1\modu p$.
Therefore $d=p-1$, so the ramification of $p$ in the singular $\Q$-field $M$ is $ e(\mk p_0/p\Z)=1$ and $e(\mk p_1/p\Z)=p-1$.
\en
%

\paragraph{Acknowledgments:}
I thank Professor Preda Mihailescu for helpful  frequent e-mail dialogues  and error detections in some intermediate versions of this paper.

%
Roland Qu\^eme

13 avenue du ch\^ateau d'eau

31490 Brax

France

mailto: roland.queme@wanadoo.fr

home page: http://roland.queme.free.fr/

************************************

MSC Classification : 11R18;  11R29; 11R32

************************************
%

\begin{thebibliography}{9}
\bibitem{hil} D. Hilbert, \textit{The Theory of Algebraic Number Fields, Zahlbericht}, Springer, 1998.
\bibitem{ire} K. Ireland, M. Rosen, \textit{A Classical Introduction to Modern Number Theory}, Springer-Verlag, 1982.
\bibitem{mol} R.A. Mollin, \textit{Algebraic Number Theory}, Chapman  and Hall/CRC, 1999.
\bibitem{nar} W. Narkiewicz, \textit{Elementary and Analytic Theory of Numbers}, Springer-verlag, 1990.
\bibitem{que} R. Qu\^eme, \textit{On $\pi$-adic expansion of singular integers of the p-cyclotomic field},
arXiv.org, preprint math.NT/0610968
\begin{verbatim}
http://arxiv.org/PS_cache/math/pdf/0610/0610968.pdf
\end{verbatim}
\bibitem{rib} P. Ribenboim, \textit{13 Lectures on Fermat's Last Theorem}, Springer-Verlag, 1979.
\bibitem{ri2} P. Ribenboim, \textit{Classical Theory of Algebraic Numbers}, Springer, 2001.
\bibitem{was} L.C. Washington, \textit{Introduction to cyclotomic fields, second edition}, Springer, 1997.
\end{thebibliography}
\end{document}